\documentclass[sn-mathphys,Numbered]{sn-jnl}

\usepackage{graphicx}%
\usepackage{multirow}%
\usepackage{amsmath,amssymb,amsfonts}%
\usepackage{amsthm}%
\usepackage{mathrsfs}%
\usepackage[title]{appendix}%
\usepackage{xcolor}%
\usepackage{textcomp}%
\usepackage{manyfoot}%
\usepackage{booktabs}%
\usepackage{algorithm}%
\usepackage{algorithmicx}%
\usepackage{algpseudocode}%
\usepackage{listings}%
\usepackage{enumerate}
\usepackage{setspace}

\theoremstyle{thmstyleone}%
\newtheorem{theorem}{Theorem}[section]%
\newtheorem{lemma}{Lemma}[section]%

\theoremstyle{thmstyletwo}%
\newtheorem{remark}{Remark}[section]%

\theoremstyle{thmstylethree}%
\newtheorem{definition}{Definition}[section]%

\numberwithin{equation}{section}
\begin{document}
	
	\title[Article Title]{Gradient Estimate for Fisher-KPP Equation on Finsler metric measure spaces}
	
	
	\author{\fnm{Bin Shen and Dingli Xia}}


\abstract{In this manuscript, we study the positive solutions of the Finslerian Fisher-KPP equation
	$$
	u_t=\Delta^{\nabla u} u+cu(1-u).
	$$
	The Fisher-KPP equation is widely applied and connected to many mathematical branches. We establish the global gradient estimates on compact Finsler metric measure manifold with the traditional $CD(K,N)$ condition, which is developed by S. Ohta and K.-T. Sturm in Finsler geometry. Furthermore, With the assistance of a new comparison theorem developed by the first author, we also give the gradient estimate on forward complete noncompact locally finite misalignment Finsler metric measure spaces with the mixed weighted curvature bounded below and some non-Riemannian curvatures norm-bounded.
	
}




\maketitle

\section{Introduction}\label{sec1}.

The Fisher-KPP equation on a complete Riemannian manifold $M$ is given by
\begin{equation}\label{1.1}
	u_{t}=\Delta  u +cu(1-u),
\end{equation}
where $u$ is a real-valued function  on $M\times[0,\infty)$ and $c$ is a positive constant.

The equation was proposed by R. A. Fisher in 1937 to describe the propagation of an evolutionarily advantageous gene in a population, and was also independently described in a seminal paper by A. N. Kolmogorov, I. G. Petrovskii, and N. S. Piskunov  in the same year; for this reason, it is often referred to in the literature as the name of Fisher–KPP equation.

Since the two papers in 1937, there have been extensive investigations on traveling wave solutions and asymptotic behavior in terms of spreading speeds for various evolution systems. Traveling waves were adopted to study the nonlinear PDEs, such as the nonlinear reaction-diffusion equations modeling physical and biological phenomena (cf. \cite{MM2002}\cite{MM2003}), the integral and integrodifferential population models (cf. \cite{AR1976}\cite{BC1977}\cite{D1978}\cite{DK1978}\cite{S1980}), the lattice differential systems (cf. \cite{BCC2003}\cite{CG2002}\cite{CG2003}\cite{CPS1998}\cite{P1999}\cite{WZ1997}), and the time-delayed reaction-diffusion equations (cf. \cite{S1987}\cite{TZ2003}).

In 2017, Cao et al. \cite{CLP2017} derived differential Harnack estimates for positive solutions to (\ref{1.1}) on Riemannian manifolds with nonnegative Ricci curvature. The idea comes from \cite{C2008}\cite{CH2009}, in which a systematic method was developed to find a Harnack inequality for geometric evolution equations. Actually, they obtained the following theorem.\\
\textbf{Theorem A.} \cite{CLP2017}
{\it   Let $(M, g)$ be an $n$-dimensional complete noncompact Riemannian manifold with nonnegative Ricci curvature, and let $u(x,t):M\times \left [ 0,\infty  \right ) \to R$ be a positive solution to (\ref{1.1}), where $u$ is $C^{2}$ in $x$ and $C^{1}$ in $t$.
	let $f=log u$, then we have 
	\begin{equation}
		\Delta f+\alpha |\nabla f|^{2}+\beta e^{f}+\phi (t)\geq 0.
	\end{equation}
	
	for all $x$ and $t$, provided that $0<\alpha<1$ as well as $\frac{-cn(2+\sqrt{2}  )}{4(1-\alpha )} <\beta<\min\{\frac{-cn(1+\alpha )}{4\alpha ^{2}-4\alpha +2n} ,\frac{-cn(2-\sqrt{2}  )}{4(1-\alpha )}\}<0$, 	
	where  	
	$\phi (t)=\frac{u(\frac{e^{2\mu wt}}{v-w}-\frac{1}{\mu+w}  )}{1-e^{2\mu wt}}, $	
	with 
	$\mu  =\beta c\sqrt{\frac{2(1-\alpha )}{c(-cn-8\beta (1-\alpha ))} },$	
	$v=\Big(\frac{4\beta (1-\alpha )}{n}+c\Big)\cdot \sqrt{\frac{2(1-\alpha )}{c(-cn-8\beta (1-\alpha ))} }$ and  
	$w=\sqrt{\frac{2(1-\alpha )}{n} }.$
}

Utilizing Theorem A, one can integrate along space-time curves to get a Harnack inequality. However, it is different from the classical Li–Yau Harnack inequality \cite{LY1986} in form.
Gradient estimates play an important role in studying elliptic and parabolic operators. The method originated first in \cite{Y1975} and \cite{CY1975}, and was further developed by Li and Yau \cite{LY1986}, Li \cite{J1991}, Negrin \cite{N1995}, P. Souplet and Q. Zhang \cite{SZ2006}, Y. Yang \cite{Y2008}, etc.. Recent gradient estimates under the geometric flow include \cite{BCP2010} and \cite{B2015}. For more results on the nonlinear PDEs, one may refer to \cite{AR2016}\cite{GR2011}.

In 2018, following the line in \cite{J1991}, X. Geng and S. Hou \cite{GH2018} proved the  positive solutions to the Fisher–KPP equation on complete Riemannian manifolds. They derived a gradient estimate, and got the classic Harnack inequality by using it, which extended the recent result of Cao, Liu, Pendleton and Ward  \cite{CLP2017}.


In 2009, upon the foundational investigations of gradient estimation on Riemannian manifolds, S. Ohta initially put forward a sophisticated framework for gradient estimation applicable to the heat equation on compact Finsler manifolds \cite{O2009}. Subsequently, in the year 2014, C. Xia \cite{X2014} undertook a comprehensive analysis of harmonic functions within the context of both compact and non-compact (specifically, forward complete) Finsler structures by utilizing an advanced form of Moser's iterative technique. Progressing further, Q. Xia  \cite{X2020} delivered intricate gradient estimates for positive solutions to the heat equation on forward complete Finsler spaces, employing methodologies analogous to those documented in \cite{X2014}. The first author \cite{S1 2023} proposed a method for global and local gradient estimates on Finsler metric measure spaces, which was also used  to solve the Finslerian Schr\"odinger equation. Later, the first author \cite{S 2023} gives the gradient estimates of bounded solutions to the Finslerian Allen-Cahn equation.

In this paper, we study the Fisher-KPP equation (\ref{1.1})
on a Finsler metric measure manifold $(M, F, \mu)$, where $c$ is a positive constant, the solution of (\ref{1.1}) is a function $u(x)$ on $M\times[0,\infty)$. Our main theorems in this paper are as follows

\begin{theorem}\label{Theorem 1.1}
	
	Let $(M, F, \mu)$ be a compact Finsler metric measure space with dimensional $n\ge2$, whose weighted Ricci curvature satisfies $Ric^N\ge-K$, for some positive constant $K$. Assume the reversibility of $M$ has upper bound $\rho_{0}$. Suppose $u$ is a bounded positive smooth solution of the  Fisher-KPP parabolic equation (\ref{1.1}) on $M\times[0,\infty)$, then we have
	
	\begin{equation}
		\frac{F^{2}(\nabla u)}{u^2}+2c(1-u) -2\frac{u_{t}}{u} \leq 2(N+n)\rho_{0}^{2}K+\frac{27}{\sqrt{2}}(N+n)^{\frac{3}{2}}M_{1}c.
	\end{equation}
	
	where, $ M_{1}=\sup u(x,t)$.
\end{theorem}

Similarily, on a noncompact forward complete Finsler manifold, the following local gradient estimates is obtained.
\begin{theorem}\label{Theorem 1.2}
	
	Let $(M,F,\mu)$ be a complete noncompact Finsler metric measure space. Denote by $B(p,2R)$ the forward geodesic ball centered at $p$ with forward radius 2$R$. Suppose the mixed weighted Ricci curvature $^{m}Ric^{N}\geq-K(2R) $  in  $B(p, 2R)$ with   $K(2R) \geq 0$, and the misalignment $\alpha$  satisfies  $\alpha \leq A(2R)$  in  $B(p, 2R)$. Moreover, the non-Riemannian tensors satisfy  $ F(U)+F^{*}(\mathcal{T} )+F(divC(V))\leq K_{0}$. Assume  $u(x,t)$ is a bounded positive smooth solution of the  Fisher-KPP parabolic equation (\ref{1.1}) on $M\times[0,\infty)$, then we have
	
	\begin{equation}
		\begin{split}
\frac{F^2_{\nabla u}(\nabla u)}{u^2}&+sc(1-u) -s\frac{u_{t}}{u}\leq C_1(N,n,\epsilon)\left(\frac{s^{2}}{t} +\frac{s^{2}K(2R)}{(s-1)} +\frac{sc}{q} M_{1}\right) \\
&+C_2(N,n,\epsilon,A(2R))\frac{s^{2}}{R^{2}} \Bigg\{1+R\left[1+\sqrt{K(2R)}\right]+\frac{s^2}{s-1}\Bigg\}.
		\end{split}
	\end{equation}
	on $B(p,R)$ for $0< \varepsilon< 1$ , $s>1$ , $q>0$, such that $\frac{2(1-\varepsilon)}{N-\varepsilon(N-n)}\frac{s-1}{sq}\geq\frac{1}{\varepsilon}-1+\frac{(2s-1)^{2}}{8}$, where $C_{1}$, $C_{2}$ are  positive constants.
	
\end{theorem}

\section{Preliminary}\label{sec2}.

In this section, we introduce the fundamental principles governing Finsler metric measure spaces (refer to \cite{S1 2023} for more details). A Finsler metric space is a triple $(M,F,\mu)$, which indicates that a differential manifold is equipped with a Finsler metric $F$ and a measure $\mu$. Suppose the local coordinate of the tangent bundle is $(x,y)$, where $x$ is the point on $M$ and $y$ is the direction on $T_{x}M$. A Finsler metric $F$ is a nonnegative function $F:TM\to \left [ 0,+\infty  \right )$ satisfies that 

(i) $F$ is smooth and positive on $TM\setminus  \left \{ 0 \right \}$;

(ii) $F$ is a positive homogenous norm, i.e., $F(x,ky)=kF(x,y)$ for any $\left ( x,y \right ) \in TM$ and for any $k>0$;

(iii) $F$ is strongly pseudo-convex, namely, for any $\left ( x,y \right ) \in TM \setminus  \left \{ 0 \right \}$, the fundamental tensor is a positive definite matrix defined by
\begin{equation}
	g_{ij}(x,y):=\frac{1}{2}\frac{\partial F^{2}}{\partial y^{i} \partial y^{j}} (x,y). 
\end{equation}

Unlike the Riemann metric, the Finsler metric is defined locally as the norm on the tangent space at each point, and globally as a metric on the pull-back bundle, so there are a large number of non-Riemannian geometric quantities on a Finsler metric measure space. The Cartan tensor is defined by
\begin{equation}
	C(X,Y,Z):=C_{ijk}X^{i}Y^{j}Z^{k}=\frac{1}{4}\frac{\partial ^{3}F^{2}(x,y)}{\partial y^{i}\partial y^{j}\partial y^{k}}X^{i}Y^{j}Z^{k},
\end{equation}
for any local vector fields $X$, $Y$, $Z$.

There is a unique almost g-compatible and torsion-free connection on the pull-back tangent bundle $\pi ^{*}TM$  of the Finsler manifold $(M, F)$ called the Chern connection. It is determined by
\begin{equation}
	\nabla _{X}Y-\nabla _{Y}X=[X,Y]; 
\end{equation}
\begin{equation}
	w(g_{y}(X,Y))-g_{y}(\nabla _{Z}X,Y) -g_{y}(X,\nabla _{Z}Y)=2C_{y}(\nabla _{Z}y,X,Y),
\end{equation}
for any $X,Y,Z \in TM \setminus  \left \{ 0 \right \}$, where $C_{y}$ is the Cartan tensor. The Chern connection coefficients are locally denoted by $\Gamma ^{i} _{jk} (x,y)$ in a natural coordinate system, which could induce the spray coefficients as $G^{i} =\frac{1}{2} \Gamma ^{i}_{jk}  y^{j} y^{k} $. The spray is given by 
\begin{equation}
	G=y^{i} \frac{\delta }{\delta x^{i}} -2G^{i}\frac{\partial }{\partial y^{i}},
\end{equation}
in which $\frac{\delta }{\delta x^{i}}=\frac{\partial }{\delta x^{i}}-N^{j}_{i}\frac{\partial}{\partial y^{j}}$  and the nonlinear connection coefficients are locally induced from the spray coefficients by $N^{i}_{j}=\frac{\partial G^{i}}{\partial y^{j}}$.

The Chern connection can define the Chern Riemannian curvature $R$ and Chern non-Riemannian connection $P$. Denote by $\Omega$ the curvature form of the Chern connection, so that
\begin{equation}
	\Omega (X,Y )Z = R(X,Y)Z + P(X,\nabla _{Y} y,Z),
\end{equation}
for any $X,Y,Z \in TM \setminus  \left \{ 0 \right \}$, where locally 
\begin{equation}
	R^{\,\,i} _{j\, kl} =\frac{\delta \Gamma ^{i}_{jl}}{\delta x^{k}}+ \frac{\delta \Gamma ^{i}_{jk}}{\delta x^{l}}+\Gamma ^{i}_{km}\Gamma ^{m}_{jl}-\Gamma ^{i}_{lm}\Gamma ^{m}_{jk} , \quad P^{\,\,i} _{j\, kl}=-\frac{\partial \Gamma ^{i}_{jk}}{\partial y^{l}}.
\end{equation}
Customarily, we denote the horizontal Chern derivative by $"\mid"$ and the vertical Chern derivative by $";"$. For example,

\begin{equation}
	v_{i\mid j}=\frac{\delta }{\delta x^{j}} v_{i}-\Gamma ^{k}_{ij}v_{k}, \quad v_{i;j}=\frac{\partial }{\partial y^{j}}v_{i},
\end{equation}
for any 1-form $v = v_{i}dx^{i}$ on the pull-back bundle.

The angular metric form $h_{y}$ is defined by
\begin{equation}
	h_{y}(u,v)=g_{y}(u,v)-\frac{1}{F^{2}(y)}g_{y}(y,u)g_{y}(y,v),
	\label{2.9}
\end{equation}
for any $y,u,v \in T_{x}M $ with $y \neq 0 $ Thus, for any two linearly independent vectors $y, v \in T_{x}M \setminus  \left \{ 0 \right \}$, which span a tangent plane $\Pi _{y}=span\{y,v\}$, the flag curvature with pole $y$ is defined by
\begin{equation}
	K(P,y)=K(y,x):=\frac{R_{y}(y,u,u,y)}{F^{2}(y)h_{y}(u,u)},
\end{equation}
which is locally expressed by
\begin{equation}
	K(y,u)=\frac{-R_{ijkl}y^{i}u^{j}y^{k}u^{l}}{(g_{ik}g_{jl}-g_{il}g_{jk}) y^{i}u^{j}y^{k}u^{l}}.
\end{equation}
The Ricci curvature is defined by
\begin{equation}
	Ric(v):=F^{2}(v)\sum_{i=1}^{n-1}K(y,e_{i}), 
\end{equation}
where $e_1, \dots, e_{n-1}, \frac{v}{F(v)}$ form an orthonormal basis of $T_x M$ with respect to $g_y$. The Landsberg curvature of $(M,F)$ is given by
\begin{equation}
	L:=L^{i}_{jk}\partial _{i}\otimes dx^{j}\otimes dx^{k},\quad L^{i}_{jk}=-y^{j}P^{\,\,i} _{j\, kl}.
\end{equation}
By the zero homogeneity, according to the Euler lemma, it easy to see that
\begin{equation}
	C_{ijk}y^{i}=L_{ijk}y^{i}=0
\end{equation}

H.B. Rademacher introduced the concepts of reversibility and reversible manifolds \cite{R2004}, which are also closely related to the analytical assumptions on Finsler manifolds. A Finsler metric is defined to be reversible if \( F(x, V) = \Bar{F}(x, V) \) for any point $x$ and any vector field $V$, where \( \Bar{F}(x, V) := F(x, -V) \) is called the reversed Finsler metric of $F$.

\begin{equation}
	\rho :=\underset{x\in M}{sup} \underset{V \neq 0}{sup} \frac{F(x,V)}{F(x,V)}.
\end{equation}

Obviously, $F$ is reversible if and only if $\rho \equiv 1$. A Finsler manifold $(M, F)$ is said to have finite reversibility if $\rho <+\infty$.

Later, K. Ball, E. Carlen, and E. Lieb introduced the concepts of uniform smoothness and uniform convexity in Banach space theory \cite{BCL2002}, whose geometric interpretation in Finsler geometry was provided by S. Ohta \cite{O2017}. We say \( F \) satisfies uniform convexity and uniform smoothness if there exist uniform positive constants \( \kappa^* \) and \( \kappa \), called the uniform convexity constant and the uniform smoothness constant, respectively, such that for any \( x \in M \), \( V \in T_x M \setminus \{ 0 \} \), and \( W \in T_x M \), we have

\begin{equation}
	\kappa ^* F^{2}(x,y)\le g_{V}(y,y)\le\kappa F^{2}(x,y).
\end{equation}

where $g_{V}=(g_{ij}(x,V))$ is the Riemannian metric on $M$ induced from $F$ with respect to the reference vector $V$ . In this situation, the reversibility $\rho$ could be controlled by $\kappa$ and \( \kappa^* \) as

\begin{equation}
	1\le \rho \le min\{\sqrt{\kappa}.\sqrt{1/\kappa^*}\}.
\end{equation}
$F$ is Riemannian if and only if $\kappa = 1$ if and only if $\kappa^* = 1$ \cite{O2017}.

The Riemannian structure is inconsistent when the reference vectors are different. For example, given three different local non-vanishing vector fields around $x$, namely, $V$, $W$, $Y$ , the norm of $Y$ about $g_{V}$ and $g_{W}$ maybe not the same in general case. The ratio $g_{V} (Y, Y )=g_{W} (Y, Y )$ is a function about $V$, $W$, $Y$ . Based on this fact, \cite{S1 2023} defined an important constant on a Finsler manifold, called the misalignment.

\begin{definition}
	(\cite{S1 2023}). For a Finsler manifold $(M,F)$, the misalignment of a Finsler metric at point $x$ is defined by
	\begin{equation}
		\alpha (x)=\underset{V,W,Y\in S_{x}M}{sup} \frac{g_{V}(Y,Y)}{g_{W}(Y,Y)} .
	\end{equation}
	Moreover, the global misalignment of the Finsler metric is defined by
	\begin{equation}
		\alpha =\underset{X\in M}{sup} \alpha_{M} (x)=\underset{X\in M}{sup}\underset{V,W,Y\in S_{x}M}{sup} \frac{g_{V}(Y,Y)}{g_{W}(Y,Y)} .
	\end{equation}
\end{definition}   

it also provided in \cite{S1 2023} some characterizations of the misalignment. Especially, a Finsler manifold $(M, F)$ is a Riemannian manifold if and only if $\alpha_{M} = 1$. Moreover, a Finsler manifold $(M,F)$ is uniform convexity and uniform smoothness if and only if it satisfies finite misalignment.

Since that, we have gaven an important class of Finsler manifold as the following.

\begin{definition}
	We call a Finsler manifold $(M,F)$ has finite misalignment if there is a positive constant $A$ such that $\alpha \le A$, and has locally finite misalignment if for any compact subset $ \Omega  \subset  M$, there is a constant $A(\Omega)$ depending on $\Omega$ such that $\alpha(x)\leq A(\omega)$ for any $x \in \Omega$.
\end{definition}

So far, we have briefly introduced some Riemannian or non-Riemannian local quantities and tensors in Finsler geometry corresponding to Riemannian geometric quantities. Next, we will introduce more tensors related to the measure $d\mu$.

For any smooth function $f:M\to R$, df  denotes its differential 1-form and its gradient $\nabla f$ is defined as the dual of the 1-form via the Legendre transformation, namely, $\nabla f(x):=l^{-1}(df(x))\in T_{x}M$. Locally it can be written as

\begin{equation}
	\nabla f =g^{ij}(x.\nabla f)\frac{\partial  f}{\partial x^{i}}\frac{\partial }{\partial x^{j}}. 
\end{equation}
on $M_{f} := {df \neq 0}$. The Hessian of $f$ is defined via the Chern connection by

\begin{equation}
	\nabla ^{2}f(X,Y)=g_{\nabla f}(\nabla _{X}^{\nabla f}\nabla f,Y).
\end{equation}
It can be shown that $\nabla^{2}f(X,Y)$ is symmetric \cite{OS2014}.

For any two points $p$, $q$ on $M$, the distance function is defined by
\begin{equation}
	d_{p}(q):=d(p,q):= \underset{\gamma}{{\inf }} \int_{0}^{1} F(\gamma (t),\dot{\gamma }(t) )dt.
\end{equation}

where the infimum is taken over all the $C^{1}$ curves $\gamma : [0, 1]\to M$ such that $\gamma(0) = p$ and $\gamma(1) = q$. Fixed a base point $p$ on $M$, we denote the forward distance function by $r$. That is, $ r(x) = d(p, x)$, with d denotes the forward distance.

The forward distance function $r$ is a function defined on the Finsler manifold $M$. $dr$ is a 1-form on $M$, whose dual is a gradient vector field, noted by $\nabla r$. Precisely, $\nabla r =g^{ij}(x,\nabla r)\frac{\partial r}{\partial x^{i}}\frac{\partial}{\partial x^{j}}$ Taking the Chern horizontal derivative of $\nabla r$ yields the Hessian of distance function $\nabla^{2}r$. Locally, in a natural coordinate system

\begin{equation}
	\nabla ^{2} r(\frac{\partial }{\partial x^{i}},\frac{\partial}{\partial x^{j} } )= \nabla _{j}\nabla _{i} r = \frac{ \partial ^{2}r}{\partial x^{i}  \partial x^{j}} -\Gamma ^{k}_{ij}(x,\nabla r)\frac{\partial r}{\partial x^{k}}.
	\label{2.23}
\end{equation}

In (\ref{2.23}), the derivative is taken in the direction $\nabla r$ naturally. Generally, we can take the derivative in any direction. Suppose $V$ is a local vector field around $x$ on $M$.

Note that the distance function may not be symmetric about $p$ and $q$ unless $F$ is reversible. $AC^{2}$ curve $\gamma$ is called a geodesic if locally

\begin{equation}
	\ddot{ \gamma} (t)+2G^{i}(\gamma(t),\dot{\gamma}(t)) =0 .
\end{equation}

where $G^{i}(x, y)$ are the spray coefficients. A forward geodesic ball centered at $p$ with radius $R$ can be represented by

\begin{equation}
	B^{+}_{R}(p):=\{q \in M:d(p,q)<R\}
\end{equation}

Adopting the exponential map, a Finsler manifold $(M,F)$ is said to be forward complete or forward geodesically complete if the exponential map is defined on the entire $TM$. Thus, any two points in a forward complete manifold $M$ can be connected by a minimal forward geodesic. Moreover, the forward closed balls $\overline{B^{+}_{R}(p)}$ are compact.

A Finsler metric measure space $(M, F, d\mu)$ is a Finsler manifold equipped with a given measure $\mu$. In local coordinates $\{x^{i}\}_{i=1}^{n}$, we can express the volume form as $d\mu = \sigma(x)dx^{1} \dots  dx^{n} $. For any  $ y \in T_{x}M \setminus  \left \{ 0 \right \}$, define

\begin{equation}
	\tau (x,y):=log \frac{\sqrt{det g_{ij}(x,y)}}{\sigma (x)} .
\end{equation}

which is called the distortion of $(M, F, d\mu)$. The definition of the S-curvature is given in the following. 

\begin{definition}\label{2.3}
	(\cite{S1997}\cite{SS2016}). Suppose $(M, F, d\mu)$ is a Finsler metric measure space. For any point $x \in M$, let $\gamma = \gamma(t)$ be a forward geodesic from x with the initial tangent vector $\dot{\gamma}(0)=y$. The S-curvature of $(M, F, d\mu)$ is 
	\begin{equation}
		S(x,y):=\frac{d}{dt}\tau =\frac{d}{dt}(\frac{1}{2} log det(g_{ij})-log \sigma (x)) (\gamma (t),\dot{\gamma}(t) )|_{t=0}.
	\end{equation}
	
	Definition \ref{2.3} means that the S-curvature is the changing of distortion along the geodesic in direction y. Modeled on the definition of $\mathcal{T}$ curvature in \cite{S2001}, \cite{S1 2023} defined the difference of $\nabla \tau$ on the tangent sphere, denoted by $\mathcal{T}$,
\end{definition}
\begin{definition}
	(\cite{S1 2023}). The difference of $\nabla \tau$ on the tangent bundle is a tensor denoted by $\mathcal{T}$, which is given by
	\begin{equation}
		\mathcal{T}(V,W):=\nabla ^{V} \tau (V)-\nabla^{W}\tau(W).
	\end{equation}
	for vector fields $V$,$W$ on $M$. Locally, it is $\mathcal{T}(V,W)=\mathcal{T}_{i}(V,w)dx^{i}$, with
	\begin{equation}
		\mathcal{T}_{i}=\frac{\delta }{\delta x^{i}} \tau _{k}(V)-\frac{\delta }{\delta x^{i}}\tau _{k}(W).
	\end{equation}
\end{definition}

Obviously, $\mathcal{T}(V, W) $ is anti-symmetric about $V$ and $W$, that is, $\mathcal{T} (V, W) =-\mathcal{T} (W, V)$ for any nonvanishing $V$ and $W$.

If in local coordinates $\{x^{i}\}^{n}_{i=1}$, expressing $d\mu= e^{\phi}dx^{1}\dots dx^{n}$, the divergence of a smooth vector field $V$ can be written as

\begin{equation}
	div_{d\mu }V=\sum_{i=1}^{n}(\frac{\partial V^{i}}{\partial x^{i}}+V^{i}\frac{\partial \phi}{\partial x^{i}}  ).
\end{equation}

The Finsler Laplacian of a function $f$ on $M$ could now be given by $\Delta _{d\mu}f:=div_{d\mu}(\nabla f)$, noticing that $\Delta_{d\mu}=\Delta ^{\nabla f}f$, where
$\Delta ^{\nabla f}_{d\mu}f:=div_{d\mu}(\nabla ^{\nabla f} f)$  is in the view of weighted Laplacian with

\begin{equation}
	\nabla ^{\nabla f}f :=\begin {cases} g^{ij}(x,\nabla f)\frac{\partial f}{ \partial x^{i}} \frac{\partial }{\partial x^{j}} \quad \quad  for \quad x\in M_{f};\\ 0  \qquad \qquad \qquad \qquad \quad for \quad x\notin M_{f} . \end{cases} 
	\end{equation}
	
	A Finsler Laplacian is better to be viewed in a weak sense due to the lack of regularity. Concretely, assuming $f \in W^{1,p}(M)$,
	
	\begin{equation}
\int_{M} \phi \Delta_{d\mu}fd\mu =-\int _{M} d \phi (\nabla f)d\mu,
\end{equation}

for any test function $\phi \in C_{0}^{\infty}(M)$.

On the other hand, the Laplacian of a function $f$ on a Riemannian manifold is the trace of the Hessian of $f$ with respect to the Riemannian metric $g$. On a Finsler metric measure space $(M,F, d\mu)$, the weighted Hessian $\tilde{H}(f)$ of a function $f$ on $M_{f} = \{x \in M : \nabla f |_{x} \neq 0\}$is defined in \cite{W2015} by

\begin{equation}
\tilde{H}(f)(X,Y)=\nabla^{2}f(X,Y)-\frac{S(\nabla f)}{n-1}h_{\nabla f}(X,Y), 
\end{equation}

where $h_{\nabla f}$ is the angular metric form in the direction $\nabla f$, given in (\ref{2.9}) It is clear that $\tilde{H}(f)$ is still a symmetric bilinear form with 

\begin{equation}
\Delta f=tr_{\nabla f}\tilde{H}(f).
\end{equation}

Inspired by this, \cite{S1 2023} defined the mixed weighted Hessian
\begin{equation}
\tilde {H}^{V}(f)(X,Y)=\nabla ^{2}(f)(X,X) -\frac{S(\nabla f)}{n-1}h_{V,\nabla f}(X,X),  
\end{equation}
where $h_{V,\nabla f}$ is the mixed angular metric form in the directions $V$ and $\nabla f$, which is defined by
\begin{equation}
h_{V,\nabla f}(X,Y)=g_{V}(X,Y)-\frac{1}{F^{2}_{V }(\nabla f)}g_{V}(X,\nabla f)g_{V}(Y,\nabla f), 
\end{equation}
for any vector $X$, $Y$.

It is necessary to remark that $ h_{\nabla f,\nabla f} = h_{\nabla f} $, so that $\tilde{H}^{\nabla f}(f)=\tilde {H}(f)$ for any function $f$ on $M$. 

With the assistance of the S-curvature, one can present the definition of the weighted Ricci curvature as the following.

\begin{definition}
(\cite{OS2014}\cite{SS2016}). Given a unit vector $V \in T_{x}M$ and an positive integral number $k$, the weighted Ricci curvature is defined by
\begin{equation}
	Ric_{k}(V):=\begin{cases}Ric(x,V)+\dot{S}(x,V)  \qquad \qquad  if\quad S(x,V)=0 \quad and\quad k=n \quad or\quad  k=\infty  \\ -\infty  \qquad \qquad \qquad \qquad \qquad \quad if \quad S(x,V)\neq 0\quad and \quad if \quad k=n\\ Ric(x,V)+\dot{S}(x,V)-\frac{S^{2}(x,V)}{k-n}  \quad if \quad n<k<\infty  \end{cases}
\end{equation}
where the derivative is taken along the geodesic started from $x$ in the direction of $V$.
\end{definition}
According to the definition of weighted Ricci curvature, B. Wu defined the weighted flag curvature when $k=N \in (1,n) \cup (n,\infty )$ in \cite{W2015}. We have completely introduced this concept for any $k$ in \cite{S1 2023}.

\begin{definition}
(\cite{S1 2023}). Let $(M, F,d\mu)$ be a Finsler metric measure space, and $V, W \in T_{x}M$ be linearly independent vectors. The weighted flag curvature $K^{k}(V , W)$ is defined by
\begin{equation}
	K^{k}(V,W):=\begin {cases} K(V;W)+\frac{\dot{S}(V) }{(n-1)F^{2}(V)}  \quad  if \quad S(x,V)=0\quad and \quad k=n \quad or\quad  k=\infty ; \\ -\infty   \qquad \qquad \qquad \qquad \quad if \quad S(x,V)\neq 0\quad and \quad if \quad k=n;\\ K(V;W)+\frac{\dot{S}(V)}{(n-1)F^{2}(V)}-\frac{S^{2}(V)}{(n-1)(k-n)F^{2}(V)}  \qquad if \quad n<k<\infty,  \end{cases}
	\end{equation}
	where the derivative is taken along the geodesic started from $x$ in the direction of $V$.
\end{definition}
Moreover, it has also been defined the mixed weighted Ricci curvature in \cite{S1 2023}, denoted by $^{m}Ric^{k}$.

\begin{definition}
(\cite{S1 2023})Given two unit vectors $V, W \in T_{x}M$ and an positive integral number $k$, the mixed weighted Ricci curvature $^{m}Ric^{k}(V,W)=^{m}Ric^{k}_{W}(V)$is defined by
\begin{equation}
^{m}Ric^{k}_{W}(V):=\begin {cases} tr_{W} R_{V}(V)+\dot{S}(x,V)\quad if \quad S(x,V)=0\quad and \quad k=n \quad or\quad  k=\infty ; \\ -\infty   \quad \qquad \qquad \qquad \quad if \quad S(x,V)\neq 0\quad and \quad if \quad k=n;\\ tr_{W}R_{V}(V)+\dot{S}(x,V)-\frac{S^{2}(V)}{k-n}  \qquad if \quad n<k<\infty , \end{cases}
\end{equation}

where the derivative is taken along the geodesic started from $x$ in the direction of $V$, and $tr_{W}R_{V}(V)=g^{ij}(W)g_{V}(R_{V}(e_{i},V)V,e_{j})$  means taking trace of the flag curvature with respect to $g(x, W)$.
\end{definition}

\begin{remark}
The weighted Ricci curvature is a special case of the mixed weighted Ricci curvature, i.e., $Ric^{k}(V)=^{m}Ric^{k}_{V}(V)$. 
Defining the function $ct_{c}(r)$ as
\begin{equation}
ct_{c}(r)=\begin{cases} \sqrt{c} cot\sqrt{c}r,  \quad \quad c>0, \\ 1/r, \qquad \qquad \quad c=0, \\\sqrt{-c}coth\sqrt{-c}r  \quad  c<0.\end{cases}
\end{equation}
the following weighted Hessian comparison theorem is cited from Theorem 3.3 in \cite{W2015}.
\end{remark}

\begin{theorem}\label{thm-2.1}
(\cite{S1 2023}). Let $(M, F, d\mu)$ be a forward complete n-dimensional Finsler metric measure space with finite misalignment $\alpha$. Denote the forward distance function by $r$ and by $V$ a fixed vector field on $M$. Suppose the mixed weighted Ricci curvature $Ric^{N} (V, \nabla r)$ of $M$ is bounded from below by $-K$ with $K > 0$, for some $N > n$, as well as the non-Riemannian curvatures $U$, $\mathcal{T}$ and $divC(V) = C^{ij}_{k|i}(V)V^{k}\frac{\delta}{\delta x^{j}}$ satisfy the norm bounds by $F(U) + F(\mathcal{T}) + F^{*}(div C(V)) \leq K_{0}$. Then, by setting $l = K/C(N, \alpha)$ with $C(N, \alpha) = N + (\alpha - 1)n - \alpha$, wherever $r$ is $C^{2}$, the nonlinear Laplacian of $r$ with reference vector $V$ satisfies

\begin{equation}
\Delta^{V}r\leq C(N,\alpha )ct_{l}(r)+\sqrt{\alpha }K_{0}.
\end{equation}
\end{theorem}

\section{Global gradient estimates on compact Finsler manifolds}\label{sec3}.

In this section, we will focus on  the global gradient estimates of positive solutions to the Fisher-KPP equation (\ref{1.1}) on compact Finsler manifolds, with the weighted Ricci curvature bounded from below. This condition of curvature is extensively employed within the realm of Finsler geometric analysis.

Suppose that $u$ is a positive solution on $M\times[0,\infty)$ to the Fisher-KPP equation(\ref{1.1}) . We define the function $W(x,t) = u^{-q}$, where  $q$  is a positive constant to be fixed later. We now have the following lemma.
\begin{lemma}
	Let $u$ be a positive solution on $M\times[0,\infty)$ to the Fisher-KPP equation (\ref{1.1}) and $W(x,t) = u^{-q}$, then $W$ satisfies that
	\begin{equation}
		(\Delta^{\nabla u}-\partial_t)W=\frac{q+1}{q}\frac{F^{2}(\nabla W)}W+cqW-cqW^{\frac{q-1}{q}}.
		\label{lemma 1}
	\end{equation}
\end{lemma}

\begin{proof}
	The gradient of $W$ with respect to $\nabla u$ is the dual of d$W$ by the pull-back metric $g_{\nabla u}$. Precisely, $\nabla W = -qu^{-q-1}\nabla u$ on $M$ and $W_{t}=-qu^{-q-1}u_{t}$, so we have
	\begin{equation}
		F^{2}_{\nabla u}(\nabla W) = q^2u^{-2q-2}F^{2}_{\nabla u}(\nabla u).
		\label{3.2}
	\end{equation}
	\begin{equation}
		tr_\nabla u (\nabla ^{\nabla u})^{2} W = q(q+1)u^{-q-2}F^{2}_{\nabla u}(\nabla u)-qu^{-q-1}tr_\nabla u (\nabla ^{\nabla u})^{2}u.
	\end{equation}
	\begin{equation}
		\frac{W_{t}}{W}=-q\frac{u_{t}}{u}.
	\end{equation}
	Then we deduce that in the distributional sense on $M_u$,
	\begin{equation}
		\begin{split}    
			\Delta ^{\nabla u}W = & div_\mu(\nabla ^{\nabla u} W) = e^{-\phi}\frac{\partial}{\partial x^{i}}(e^\phi g^{ij}W_i)\\ 
			= & tr_\nabla u (\nabla ^{\nabla u})^{2}W-S(\nabla W) \\
			=& q(q+1)u^{-q-2}F^{2}_{\nabla u}(\nabla u)-qu^{-q-1}\Delta^{\nabla u} u.
		\end{split}
		\label{Delta W}
	\end{equation}
	wherever $\nabla u \neq 0$ Combine with (\ref{1.1}) and (\ref{3.2}), we have
	\begin{equation}
		\begin{split}
			\Delta^{\nabla u}W= &\frac{q+1}{q}\frac{F^{2}_{\nabla u}(\nabla W)}{W}-qu^{-q-1}[u_{t}-cu(1-u)]\\
			=&\frac{q+1}{q}\frac{F^{2}_{\nabla u}(\nabla W)}{W}+W_{t}+cqW-cqW^{\frac{q-1}{q}}.
		\end{split}
	\end{equation}
	
	After the transformation, we obtain the Lemma \ref{lemma 1}.\\
\end{proof}

Now, we follow the line in \cite{J1991}. Define three functions as the follow.
\begin{equation}
	\begin{cases}
		H_{0}(x,t)=\frac{F^{2}_{\nabla u}(\nabla W)}{W^{2}}+ac(1-W^{-\frac{1}{q}}),\\
		H_{1}=\frac{W_{t}}{W},\\
		H=H_{0}+\beta H_1,\\
	\end{cases}
\end{equation}
where $a$, $\beta$ are two positive constants to be determined later. Direct calculations provide that
\begin{equation}
	\nabla^{\nabla u}H_0=\frac{\nabla^{\nabla u}F^{2}_{\nabla u}(\nabla^{\nabla u} W)}{W^2}-\frac{2F^{2}_{\nabla u}(\nabla ^{\nabla u}W)\nabla ^{{\nabla u}} W}{W^{3}}+\frac{ac}{q}W^{-\frac{q+1}{q}}\nabla ^{\nabla u}W,
\end{equation}
and
\begin{equation}
	\begin{split} 
		tr_{\nabla u}(\nabla ^{\nabla u})^{2} H_0 =& \frac{2}{W^{2}} [dW(\nabla^{\nabla u}(tr_{\nabla u}(\nabla ^{\nabla u})^{2} W)) + Ric(\nabla ^{\nabla u }W) + \left\| (\nabla^{\nabla u})^{2} W \right\|^{2}_{HS(\nabla u)}]\\
		&-\frac{8}{W^{3}} g_{\nabla u} ((\nabla ^{\nabla u})^{2} W,\nabla ^{\nabla u}W\otimes\nabla ^{\nabla u} W) 
		+ \frac{6F^{4}_{\nabla u}(\nabla ^{\nabla u}W)}{W^{4}} \\
		&-\frac{2F^{2}_{\nabla u}(\nabla ^{\nabla u} W)tr_{\nabla u} (\nabla^{\nabla u}) ^{2}W}{W^{3}} - \frac{ac(q+1)}{q^2}\frac{W^{-\frac{1}{q}}}{W^{2}}F^{2}_{\nabla u }(\nabla ^{\nabla u} W) \\
		&+\frac{ac}{q}W^{-\frac{q+1}{q}}tr_{\nabla u} (\nabla^{\nabla u}) ^{2}W,
		\label{tr_H_0}
	\end{split}
\end{equation}
on $M_{u}$, where we have employed the Ricci-type identity in \cite{S2018}. On the other hand, it satisfies that
\begin{equation}
	\nabla^{\nabla u}(S(\nabla ^{\nabla u} W))=\dot{S}(\nabla ^{\nabla u} W)+g_{\nabla u}((\nabla^{\nabla u})^{2}W,\nabla^{\nabla u} \tau\otimes\nabla ^{\nabla u} W),
	\label{3.10}
\end{equation}
and
\begin{equation}
	\nabla^{\nabla u}(tr_{\nabla u}(\nabla ^{\nabla u})^{2}W)=\nabla^{\nabla u}(\Delta^{\nabla u} W)+\dot{S}(\nabla ^{\nabla u}W)+g_{\nabla u}((\nabla^{\nabla u})^{2}W,\nabla^{\nabla u} \tau\otimes\nabla ^{\nabla u} W) .
	\label{3.11}
\end{equation}
Plugging (\ref{3.10}) and (\ref{3.11}) into (\ref{tr_H_0}) yields
\begin{equation}
	\begin{split} 
		tr_{\nabla u}(\nabla ^{\nabla u})^{2} H_0=&\frac{2}{W^{2}} [dW(\nabla^{\nabla u}(\Delta ^{\nabla u} W) +g_{\nabla u}((\nabla^{\nabla u})^{2}W,\nabla^{\nabla u} \tau\otimes\nabla ^{\nabla u} W)\\
		 &+\dot{S}(\nabla ^{\nabla u} W)+ Ric(\nabla ^{\nabla u} W)
		+ \left\| (\nabla^{\nabla u})^{2} W \right\|^{2}_{HS(\nabla u)}]\\
		 &+\frac{6 F^{4}_{\nabla u}(\nabla ^{\nabla u}W)}{W^{4} }- \frac{2F^{2}_{\nabla u}(\nabla ^{\nabla u}W)(\Delta ^{\nabla u} W+S(\nabla ^{\nabla u}W))}{W^{3}}  \\
		 &-\frac{8}{W^{3}} g_{\nabla u} ((\nabla ^{\nabla u})^{2}W,\nabla ^{\nabla u}W\otimes \nabla ^{\nabla u}W)  \\
		 &-\frac{ac(q+1)}{q^2}\frac{W^{-\frac{1}{q}}}{W^{2}}F^{2}_{\nabla u}(\nabla ^{\nabla u} W) +\frac{ac}{q}W^{-\frac{q+1}{q}}(\Delta ^{\nabla u} W+S(\nabla ^{\nabla u} W)),
		\label{trH_01}
	\end{split}
\end{equation}
wherever $\nabla u \neq 0$.
By a basic computation of Finsler geometry, we could get
\begin{equation}
	\begin{split}
		\Delta ^{\nabla u}H_0=div_\mu(\nabla^{\nabla}H_0)= tr_{\nabla u}(\nabla ^{\nabla u})^{2}H_0 -d\tau(\nabla^{\nabla u}H_0)+2C_{\nabla^{2}u}^{\nabla u}(\nabla^{\nabla u}H_0),
		\label{H_01}
	\end{split}
\end{equation}
where we have already employed the fact that $C(\nabla u,\cdot,\cdot)=0 $ and $C_{\nabla^{2}u}^{\nabla u}(\nabla^{\nabla u}H_0)=u_{|i}^{k}C_{k}^{ij}(\nabla u)H_j$. The notation ``$|$'' denotes the horizontal derivative with respect to the Chern connection in the direction $\nabla u$. It follows from
\begin{equation}
	\begin{split} 
		C_{\nabla^{2}u}^{\nabla u}(\nabla^{\nabla u}H_0)=&u_{|i}^{k}C_{k}^{ij}(\nabla u)[2W^{l}W_{l|j}-2F^{2}(\nabla W)W^{-3}W_j+\frac{1}{q}acW^{-\frac{q+1}{q}}W_j] \\
		=&\frac{2u_{|i}^{k}}{W^2}[(C_{k}^{ij}W_j)_{|l}W^{l}-C_{k|l}^{ij}W_jW^{l}]\\
		=&-2\frac{u_{|i}^{k}}{W^{2}}L^{ij}_{k}W_j\\=&0,
	\end{split}
\end{equation}
that (\ref{H_01}) is equal to  
\begin{equation}
	\Delta ^{\nabla u}H_0=div_\mu(\nabla^{\nabla u}H_0)= tr_{\nabla u} (\nabla^{\nabla u})^2H_0 -d\tau(\nabla^{\nabla u}H_0).     
\end{equation}

Combining it with (\ref{trH_01}), we could find that 
\begin{equation}
	\begin{split}
		\Delta^{\nabla u}H_0&=\frac{2}{W^{2}} [g_{\nabla u}(\nabla ^{\nabla u}(\Delta ^{\nabla u}W),\nabla ^{\nabla u}W)+ Ric^{\infty}(\nabla W)  + \left\| (\nabla^{\nabla u})^{2} W \right\|^{2}_{HS(\nabla u)}]\\&-\frac{8}{W^{3}} g_{\nabla u} ((\nabla ^{\nabla u})^{2}W,\nabla^{\nabla u} W\otimes\nabla^{\nabla u} W) 
		- \frac{2F^{2}_{\nabla u}(\nabla^{\nabla u} W)(\Delta ^{\nabla u}W)}{W^{3}} \\
		&+ \frac{6F^{4}_{\nabla u}(\nabla ^{\nabla u} W)}{W^{4}} - \frac{ac(q+1)}{q^2}\frac{W^{-\frac{1}{q}}}{W^{2}}F^{2}_{\nabla u}(\nabla ^{\nabla u}W) +\frac{ac}{q}W^{-\frac{q+1}{q}}\Delta^{\nabla u} W,
		\label{H001}
	\end{split}
\end{equation}
by noticing that 
$$d\tau(\nabla^{\nabla u}H)=2d\tau[\frac{g_{\nabla u}((\nabla^{\nabla u})^{2}W,\nabla^{\nabla u}W)}{W^{2}}-\frac{F^{2}_{\nabla u}(\nabla ^{\nabla u}W)}{W^{3}}\nabla^{\nabla u }W+\frac{ac}{2q}W^{-\frac{q+1}{q}}\nabla ^{\nabla u}W].$$

Therefore, $ H_{0}(x,t)$ satisfies that
\begin{equation}
	\frac{\partial H_{0}(x,t)}{\partial_t}=\frac{\partial_t F^{2}_{\nabla u}(\nabla W)}{W^2} -\frac{2 F^{2}_{\nabla u}(\nabla^{\nabla u} W)W_t}{W^3}+\frac{ac}{q}W^{-\frac{q+1}{q}}W_{t}.
	\label{H002}
\end{equation}
Because of
$\nabla W_{t}=(\nabla W)_{t}$,
it follows for $\nabla^{\nabla u} H_1$ that
\begin{equation}
	\nabla^{\nabla u} H_1 = \nabla(\frac{W_t}{W})=\frac{\nabla W_t}{W} - \frac{W_t\nabla W}{W^2}.
\end{equation}
It asserts by direct caiculation that
\begin{equation}
	\frac{\partial }{\partial_t} (\Delta ^{\nabla u} W)=\Delta ^{\nabla u} W_{t},
\end{equation}
arccording to \cite{S1 2023}. Thus,
for $H_1$, one may find that
\begin{equation}
	\Delta ^{\nabla u} H_1 = \frac{W(\Delta  ^{\nabla u} W)_t-2(\nabla ^{\nabla u} W)_t\nabla ^{\nabla u} W-\Delta ^{\nabla u}W W_t}{W^2}+\frac{2W_tF^{2}_{\nabla u}(\nabla ^{\nabla u} W)}{W^3},
\end{equation}
so that
\begin{equation}
	\begin{split}
		\left(\Delta ^{\nabla u} -\partial_t\right)H_1
		&=\frac{2}{q}g_{\nabla u}(\nabla ^{\nabla u} H_{1},\nabla ^{\nabla u} \log W )+ cW^{-\frac{1}{q}-1}W_t.
	\end{split}
	\label{3.22}
\end{equation}
Moreover, one could calculate that
\begin{equation}
	\begin{split}
		&\frac{2}{W^2}g_{\nabla u}(\nabla ^{\nabla u}(\Delta ^{\nabla u}W),\nabla ^{\nabla u}W)-\frac{1}{W^2}\partial_t F^2(\nabla W)\\
		=&\frac{4(q+1)}{q}\frac{g_{\nabla u}((\nabla^{\nabla u})^{2}W,\nabla^{\nabla u} W\otimes\nabla^{\nabla u} W)}{W^3}\\
		&+\frac{2cF^{2}_{\nabla u}(\nabla^{\nabla u} W)}{W^2}[q-(q-1)W^{-\frac{1}{q}}]-\frac{2(q+1)}{q}\frac{F^{4}_{\nabla u}(\nabla ^{\nabla u}W)}{W^4},
		\label{cha 1}
	\end{split}
\end{equation}
\begin{equation}
	\begin{split}
		\frac{-2F^{2}_{\nabla u}(\nabla ^{\nabla u}W)}{W^3}(\Delta ^{\nabla u}-\partial_t)W 
		= &-\frac{2(q+1)}{q}\frac{F^{4}_{\nabla u}(\nabla ^{\nabla u}W)}{W^4}-2cq\frac{F^{2}_{\nabla u}(\nabla ^{\nabla u}W)}{W^2}\\
		&+2cqW^{-\frac{1}{q}}\frac{F^2_{\nabla u}(\nabla ^{\nabla u}W)}{W^2},
		\label{cha 2}
	\end{split}
\end{equation}
as well as
\begin{equation}
	\begin{split}
		\frac{ac}{q}W^{-\frac{q+1}{q}}(\Delta ^{\nabla u}-\partial_t)W =\frac{ac(q+1)}{q^2}\frac{F^{2}_{\nabla u}(\nabla ^{\nabla u}W)}{W^2}W^{-\frac{1}{q}}+ac^2W^{-\frac{1}{q}}-ac^2W^{-\frac{2}{q}},
		\label{cha 3}
	\end{split}
\end{equation}

Plugging (\ref{cha 1})-(\ref{cha 3}) into (\ref{H001}) and (\ref{H002}) yields
\begin{equation}
	\begin{split}
		\left(\Delta^{\nabla u} -\partial_t\right)H_{0} =&
		\frac{4(1-q)}{q}\frac{g_{\nabla u}((\nabla ^{\nabla u})^{2}W,\nabla^{\nabla u} W\otimes\nabla ^{\nabla u}W)}{W^3}+2\frac{q-2}{q}\frac{F^{4}_{\nabla u}(\nabla^{\nabla u} W)}{W^4}\\
		&+2 \frac{\|(\nabla^{\nabla u})^{2}W\|^{2} _{HS(\nabla u)}}{W^2}+ac^{2}W^{-\frac{1}{q}}(1-W^{-\frac{1}{q}})+\frac{2Ric^{\infty}(\nabla ^{\nabla u }W)}{W^2}\\
		&+ 2cW^{-\frac{1}{q}}\frac{F^{2}_{\nabla u}(\nabla^{\nabla u} W)}{W^2}.
	\end{split}
	\label{3.26}
\end{equation}
It asserts from the h\"{o}lder inequality that
\begin{equation}
	\frac{2g_{\nabla u}((\nabla ^{\nabla u})^{2}W,\nabla^{\nabla u} W\otimes\nabla ^{\nabla u}W)}{W^3} \leq \frac{\varepsilon \|(\nabla^{\nabla u})^{2}W\|^{2} _{HS(\nabla u)}}{W^{2}}+\frac{1}{\varepsilon}\frac{F^{4}_{\nabla u}(\nabla ^{\nabla u} W)}{W^{4}}.
\end{equation}
Thus, when $0<\varepsilon<1$, it satisfies that
\begin{equation}
	\begin{split}
		&\frac{4(1-q)}{q}\frac{g_{\nabla u}((\nabla ^{\nabla u})^{2}W,\nabla^{\nabla u} W\otimes\nabla ^{\nabla u}W)}{W^{3}}+\frac{2\left\|(\nabla^{\nabla u})^{2}W\right\|^{2}_{HS(\nabla u)}}{W^{2}}+\frac{2(q-2)}{q}\frac{F^{4}_{\nabla u}(\nabla W)}{W^{4}}\\ 
		\geq&\frac{4}{q}\frac{g_{\nabla u}((\nabla ^{\nabla u})^{2}W,\nabla^{\nabla u} W\otimes\nabla ^{\nabla u}W)}{W^{3}}+2(1-\varepsilon)\frac{\left\|(\nabla^{\nabla u})^{2}W\right\|^{2}_{HS(\nabla u)}}{W^{2}}\\&+2(1-\frac{1}{\varepsilon}-\frac{2}{q})\frac{F^{4}_{\nabla u}(\nabla ^{\nabla u}W)}{W^{4} }. 
		\label{3.28}
	\end{split}
\end{equation}

According to the equality $\frac{(a+b)^{2} }{n} =\frac{a^{2}}{N} -\frac{b^{2}}{N-n}+\frac{N(N-n)}{n}(\frac{a}{N}+\frac{b}{N-b} )^{2}$, for any $N > n$, the following inequality holds by substituting $a$, $b$ and $N$ by $\Delta ^{\nabla u} W$, $\dot{S}(\nabla ^{\nabla u}W)$ and $ N -\varepsilon(N -n)$, respectively. Namely,
\begin{equation}
	\begin{split}
		\left\|(\nabla^{\nabla u})^{2}W\right\|^{2}_{HS(\nabla u)} &\geq \frac{(tr_{\nabla u}(\nabla ^{\nabla u} )^{2} W)^{2}}{n} =\frac{(\Delta ^{\nabla u} W+S(\nabla ^{\nabla u}W) )^{2} }{n} \\
		&\geq\frac{(\Delta^{\nabla u} W)^{2} }{N-\varepsilon (N-n)} -\frac{(S(\nabla ^{\nabla u}W))^{2} }{(1-\varepsilon )(N-n)},
	\end{split}
\end{equation}
so that
\begin{equation}
	\begin{split}
		&2(1-\varepsilon)\frac{\left\|(\nabla^{\nabla u})^{2}W\right\|^{2}_{HS(\nabla u)} }{W^{2}}+\frac{2Ric^{\infty}(\nabla ^{\nabla u}W)}{W^{2} }\\
		\geq&\frac{2(1-\varepsilon )}{N-\varepsilon (N-n)}\frac{(\Delta^{\nabla u}  W)^{2}}{W^{2} }-\frac{2S((\nabla^{\nabla u} W))^2}{(N-n)W^{2} }+ \frac{2Ric^{\infty}(\nabla ^{\nabla u}W)}{W^{2} }\\
		\geq&\frac{2(1-\varepsilon )}{N-\varepsilon (N-n)}\frac{(\Delta^{\nabla u}  W)^{2}}{W^{2} }+\frac{2Ric^{N}(\nabla ^{\nabla u}W)}{W^{2} }.    	
	\end{split}
\end{equation}
Moreover, it's easy to infer that
\begin{equation}
	\begin{split}
		g_{\nabla u}  (\nabla ^{\nabla u }  H_{0},\nabla ^{\nabla u} \log W)=&\frac{2g_{\nabla u}((\nabla ^{\nabla u})^{2}W,\nabla ^{\nabla u}W\otimes \nabla ^{\nabla u}W )}{W^{3} }-\frac{2F^{4}_{\nabla u}(\nabla^{\nabla u} W)}{W^{4} }\\
		&+\frac{acW^{-\frac{1}{q} } }{q} \frac{F^{2}_{\nabla u}(\nabla^{\nabla u} W)}{W^{2} }.
	\end{split}
	\label{3.31}
\end{equation}

It is deduced from (\ref{3.26}) by combining (\ref{3.28})-(\ref{3.31}) and employing the definition  of weighted Ricci curvature that

\begin{equation}
	\begin{split}
		\left(\Delta ^{\nabla u}-\partial_t\right)H_{0} \geq& \frac{2(1-\varepsilon)}{N-\varepsilon(N-n)}\frac{(\Delta^{\nabla u}W)^2}{W^{2}}+\frac{2Ric^N(\nabla ^{\nabla u} W)}{W^2}\\&+\frac{2}{q}g_{\nabla u}\left(\nabla^{\nabla u}H_0,\nabla ^{\nabla u}log W\right)
		+2\left(1-\frac{1}{\varepsilon }\right)\frac{F^{4}_{\nabla u}\left(\nabla ^{\nabla u}W\right)}{W^{4}}\\&+ac^2
		W^{-\frac{1}{q}}(1-W^{-\frac{1}{q}})+2c\frac{F^{2}_{\nabla u}(\nabla^{\nabla u} W)}{W^2}W^{-\frac{1}{q}}\left(1-\frac{a}{q^2}\right),
	\end{split}
\end{equation}
on $M_{u}$.

Let $\beta = \frac{a}{q}$, combining with (\ref{3.22}), we have
\begin{equation}
	\begin{split}
		\left(\Delta^{\nabla u} -\partial_t\right)H \geq& \frac{2(1-\varepsilon)}{N-\varepsilon(N-n)}\frac{(\Delta^{\nabla u}W)^2}{W^2} -2\left(\frac{1}{\varepsilon}-1\right)\frac{F^{4}_{\nabla u}\left(\nabla^{\nabla u} W\right)}{W^4}\\&+\frac{2}{q}g_{\nabla u}\left(\nabla^{\nabla u}H,\nabla ^{\nabla u}\log W\right)+2c\left(1-\frac{a}{q^2}\right)\frac{F^{2}_{\nabla u}(\nabla^{\nabla u} W)}{W^2}W^{-\frac{1}{q}}\\&+ac^2W^{-\frac{1}{q}}(1-W^{-\frac{1}{q}})+\frac{ac}{q}W^{\frac{1}{q}-1}W_t+\frac{2Ric^N(\nabla ^{\nabla u}W)}{W^2}.
	\end{split}
	\label{3.33}
\end{equation}
Noticing $\frac{\Delta ^{\nabla u} W}{W} =\frac{q}{a} H+({\frac{q+1}{q} }-\frac{q}{a}  )\frac{F^{2}_{\nabla u}(\nabla ^{\nabla u}W) }{  W^{2} }$, and setting $a=sq^{2}$, we get that
\begin{equation}
	\begin{split}
		\frac{\Delta ^{\nabla u}W}{W} = \frac{1}{sq}H+\left(\frac{q+1}{q}-\frac{1}{sq}\right)\frac{F^{2}_{\nabla u}(\nabla^{\nabla u} W)}{W^2} = \frac{1}{sq}H+\left(\frac{q+1-\frac{1}{s}}{q}\right)\frac{F^{2}_{\nabla u}(\nabla ^{\nabla u}W)}{W^2}.
	\end{split}
\end{equation}
Plugging it into (\ref{3.33}), we arrive at the following lemma.

\begin{lemma}\label{lem-1}
	Let $(M, F, \mu)$ be a  forward complete Finsler metric measure space, and denote $M_{u}=\{x \in M \mid \nabla u(x) \neq 0\}$. For $H=\frac{F^{2}_{\nabla u}(\nabla ^{\nabla u}W)}{W^{2}}+ac(1-W^{-\frac{1}{q}})+\beta \frac{W_{t}}{W}$ with $a=sq^{2}$, $\beta=\frac{a}{q}$, it satisfies on $M_{u}$ that
	\begin{equation}
		\begin{split}
			\left(\Delta^{\nabla u}-\partial_t\right)H \geq& \left[\frac{2(1-\varepsilon)}{N-\varepsilon(N-n)}\frac{(sq+s-1)^{2}}{s^2q^2}-2\left(\frac{1}{\varepsilon}-1\right)\right]\frac{F^{4}_{\nabla u}(\nabla ^{\nabla u}W)}{W^4}\\
			&+ \frac{2(1-\varepsilon)}{N- 
				\varepsilon(N-n)}\frac{1}{s^2q^2}H^2 +\frac{4(1-\varepsilon)}{N-\varepsilon(N-n)}\frac{sq+s-1}{s^2q^2}H\frac{F^{2}_{\nabla u}\left(\nabla^{\nabla u} W\right)}{W^2}\\&+\frac{2}{q}g_{\nabla u}\left(\nabla^{\nabla u}H,\nabla ^{\nabla u}\log W\right)+\left(c-2cs\right)\frac{F^{2}_{\nabla u}(\nabla ^{\nabla u}W)}{W^2}W^{-\frac{1}{q}}\\
				&+cW^{-\frac{1}{q}}H+\frac{2Ric^N(\nabla ^{\nabla u}W)}{W^2} .
		\end{split}
		\label{lemma 2}
	\end{equation}
	\label{lemma}
\end{lemma}

Employing Lemma \ref{lem-1}, the following global gradient estimate theorem could be obtained by utilizing the maximum principle argument, which was first adopted in \cite{OS2014} and was applied also in \cite{S1 2023}\cite{S 2023}.

\begin{theorem}
	Let $(M, F, \mu)$ be a compact Finsler metric measure manifold  whose weighted Ricci curvature satisfies $Ric^N\ge-K$, for some positive constant $K$.  Assume the bound of the reversibility on $M$ is $\rho_{0}$. Suppose $u(x,t)$ is a bounded positive smooth solution of the  Fisher-KPP parabolic equation (\ref{1.1}) on $M\times[0,\infty)$, then we have	
	\begin{equation}
		\begin{split}
			\frac{F^2(\nabla u)}{u^2}+sc(1-u) -s\frac{u_{t}}{u}
			\leq\frac{[N-\varepsilon(N-n)]}{2(1-\varepsilon)(s-1)}s^{2}K+\sqrt{\frac{N-\varepsilon(N-n)}{2(1-\varepsilon)}}\frac{s}{q}M_1c.
		\end{split}
	\end{equation}	
	where, $0< \varepsilon< 1$, $s>1$, $q>0$, such that $\frac{2(1-\varepsilon)}{N-\varepsilon(N-n)}\frac{s-1}{sq}\geq\frac{1}{\varepsilon}-1+\frac{(2s-1)^{2}}{8}$, and $ M_{1}=\sup u(x,t)$.

	\label{theorem of compact}
\end{theorem}

\begin{proof}

Lemma \ref{lemma} implies that (\ref{lemma 2}) holds in the distributional sense on $M$. For any nonnegative test function $\varphi$, we have

\begin{equation}
	-\int_{0}^{T} \int_{M}\left [ d \varphi (\nabla H)- \varphi_{t}H \right ]d \mu dt \geq \int_{0}^{T} \int _{M} \varphi \beta d\mu dt.
	\label{ruo 0}
\end{equation}
where $\beta$ denotes the RHS of (\ref{lemma 2}). 

Assume $y_0=(x_0,t_0)\in B(p,2R)\times[0,\infty)$ be the point where $H$ achieves its maximum. Without loss of generality, we could assume $H(x_{0},t_{0}) \geq 0$, otherwise the result will be satisfied trivially. We claim that $\beta(x_{0},t_{0}) \leq 0$. 

Otherwise, $\beta$ is strictly positive at $(x_{0},t_{0})$, so that $\beta(x_{0},t_{0})$ is positive in a small neighborhood of $(x_{0},t_{0})$ on $M$, which may be denoted by $U$. Chosen a test function $\varphi$ whose compact support $V$ is contained in $U$, we know from (\ref{ruo 0}) that $H$ is a weak, local subharmonic function in a neighborhood $V\subset U$. It is a contradiction because $(x_{0},t_{0})$ is a inner point of $V$.

Because of $\beta(x_{0},t_{0})\leq 0$ and $\nabla ^{\nabla u}H(y_{0})=0$, we arrive at
\begin{equation}
	\begin{split}
		0\geq&\frac{2(1-\varepsilon)}{N-\varepsilon(N-n)}\frac{1}{s^{2}q^{2}}H^{2}+\left[\frac{2(1-\varepsilon)}{N-\varepsilon(N-n)}\frac{(sq+s-1)^2}{s^{2}q^{2}} - 2(\frac{1}{\varepsilon}-1)\right]\frac{F^{4}_{\nabla u}(\nabla ^{\nabla u}W)}{W^{4}}\\
		&+\frac{4(1-\varepsilon)}{N-\varepsilon(N-n)}\frac{sq+s-1}{s^{2}q^{2}}H\frac{F^{2}_{\nabla u}(\nabla^{\nabla u} W)}{W^{2}}+(c-2cs)\frac{F^{2}_{\nabla u}(\nabla ^{\nabla u}W)}{W^{2}}W^{-\frac{1}{q}}\\
		&+cW^{-\frac{1}{q}}H+\frac{2Ric^N(\nabla ^{\nabla u}W)}{W^{2}}.
	\end{split}
	\label{3.38}
\end{equation}
It follows from  the h\"{o}lder inequality that 
\begin{equation}
	\begin{split}
		\frac{2Ric^N(\nabla ^{\nabla u}W)}{W^{2}}&\geq-\frac{2K\rho _{0}^{2}F_{\nabla u}^{2}(\nabla ^{\nabla u}W)}{W^{2}}\\
		&\geq-\left[\frac{2(1-\varepsilon)(s-1)^{2}}{N-\varepsilon(N-n)}\frac{F_{\nabla u}^{4}(\nabla ^{\nabla u}W)}{s^{2}q^{2}W^{4}}+\frac{(N-\varepsilon(N-n))s^{2}q^{2}}{2(1-\varepsilon)(s-1)^{2}}\rho_{0}^{4}K^{2}\right],
	\end{split}
\end{equation}
and
\begin{equation}
	\begin{split}
		(c-2cs)\frac{F^{2}_{\nabla u}(\nabla ^{\nabla u}W)}{W^{2}}W^{-\frac{1}{q}}&\geq -\left[\frac{(2s-1)^2}{4}\frac{F_{\nabla u}^{4}(\nabla ^{\nabla u}W)}{W^{4}}+W^{-\frac{2}{q}}c^{2}\right]\\&\geq-\left[\frac{(2s-1)^{2}}{4}\frac{F^{4}_{\nabla u}(\nabla ^{\nabla u}W)}{W^{4}}+M_{1}^{2}c^{2}\right],
	\end{split}
\end{equation}
where $M_{1}=\sup u(x)$. Now, let $q>0$ such that $\frac{2(1-\varepsilon)}{N-\varepsilon(N-n)}\frac{s-1}{sq}\geq\frac{1}{\varepsilon}-1+\frac{(2s-1)^{2}}{8}$. Combining with (\ref{3.38}), we get 
\begin{equation}
	\begin{split}
		\frac{2(1-\varepsilon)}{N-\varepsilon(N-n)}\frac{1}{s^{2}q^{2}}H^{2}&\leq \frac{(N-\varepsilon(N-n))s^{2}q^{2}}{2(1-\varepsilon)(s-1)^{2}}\rho_{0}^{4}K^{2}+M_{1}^{2}c^{2}\\
		&\quad-\frac{4(1-\varepsilon)(sq+s-1)}{[N-\varepsilon(N-n)]s^2q^{2}}H\frac{F_{\nabla u}^{2}(\nabla ^{\nabla u}W)}{W^{2}}-cW^{-\frac{1}{q}}H\\
		&\leq \frac{(N-\varepsilon(N-n))s^{2}q^{2}}{2(1-\varepsilon)(s-1)^{2}}\rho_{0}^{4}K^{2}+ M_{1}^{2}c^{2}.
	\end{split}
\end{equation}

Hence
\begin{equation}
	\begin{split}
		H\leq\frac{[N-\varepsilon(N-n)]}{2(1-\varepsilon)(s-1)}s^{2}q^{2}\rho_{0}^{2}K+\sqrt{\frac{N-\varepsilon(N-n)}{2(1-\varepsilon)}}M_1csq.
	\end{split}
\end{equation}
we achieve the global gradient estimates on compact Finsler manifolds. 
\end{proof}

\begin{remark}
	Theorem  \ref{Theorem 1.1} follows by taking  $\varepsilon=\frac{1}{2}$ , $s=2$, $q=\frac{2}{27(N+n)}$ in Theorem \ref{theorem of compact} .
\end{remark}

\section{Local gradient estimates on forward complete Finsler spaces}

 In this section, we prove the local gradient estimates on forward complete Finsler metric measure spaces with the assistance of Lemma \ref{lemma} and the Comparison theorem (cf. Theorem \ref{thm-2.1}).

\begin{theorem}[]
	Let $(M,F,\mu)$ be a complete noncompact Finsler metric measure space. Denote by $B(p,2R)$ the forward geodesic ball centered at $p$ with forward radius $2R$. Suppose the mixed weighted Ricci curvature $^{m}Ric^{N}\geq-K(2R) $  in  $B(p, 2R)$ with  $K(2R) \geq 0$, and the misalignment $\alpha$  satisfies  $\alpha \leq A(2R)$  in  $B(p, 2R)$. Moreover, the non-Riemannian tensors satisfy  $ F(U)+F^{*}(\mathcal{T} )+F(divC(V))\leq K_{0}$. Assume  $u(x,t)$ is a bounded positive smooth solution of the  Fisher-KPP parabolic equation (\ref{1.1}) on $M\times[0,\infty)$, then we have
	
	\begin{equation}
		\begin{split}
			&\frac{F^2_{\nabla u}(\nabla u)}{u^2}+sc(1-u) -s\frac{u_{t}}{u}\\
			\leq &\frac{N-\varepsilon(N-n)}{2(1-\varepsilon)}\frac{s^{2}}{t} +\frac{N-\varepsilon (N-n)}{2(1-\varepsilon )(s-1)} s^{2}K(2R)+\frac{sc}{q} \sqrt{\frac{N-\varepsilon (N-n)}{2(1-\varepsilon )}}M_{1} \\
			&+\frac{[N-\varepsilon (N-n)]s^{2}}{2(1-\varepsilon )R^{2}} \Bigg\{C_3[1+R(1+\sqrt{K(2R)})]+\frac{[N-\varepsilon(N-n)]s^2}{4(1-\varepsilon)(s-1)}C^{2}_{1} +2C_{1}^2\Bigg\}. 			
		\end{split}
	\end{equation}
	on $B(p,R)$ for $0< \varepsilon< 1$, $s>1$, $q>0$, such that $\frac{2(1-\varepsilon)}{N-\varepsilon(N-n)}\frac{s-1}{sq}\geq\frac{1}{\varepsilon}-1+\frac{(2s-1)^{2}}{8}$, where $C_{1}$, $C_{2}$ are positive constants.
	
	\label{theorem of noncompact}
\end{theorem}

\begin{proof}

In noncompact situation, we choose the cut-off function $\tilde{\varphi}\in C^{2}[0,+\infty)$, such that
\begin{equation}
	\begin{cases}\tilde{\varphi}(r)=1&r\leq1\\\tilde{\varphi}(r)=0&r\geq2.\end{cases}   
\end{equation}
So $\tilde\varphi(r)\in[0,1]\nonumber$. 
In addition, we supposed that
\begin{equation}
	-C_{1}\leq\frac{\tilde{\varphi}^{\prime}(r)}{\tilde{\varphi}^{\frac{1}{2}}(r)}\leq0 , \quad -C_{2}\leq \tilde{\varphi}^{\prime\prime}(r),
\end{equation}
where $C_{1}, C_{2}$ are positive constants.

For a fixed point $p$, denote by $r(x)$ the forward distance function from $p$ to any point $x$. We define the cut-off function by $\varphi(x)=\tilde\varphi(\frac{r(x)}{R}) $.
So that
\begin{equation}
	F_{\nabla u}(\nabla ^{\nabla u}\varphi)\leq \frac{\sqrt{\alpha}C_1}{R}\varphi^{\frac{1}{2}}.
\end{equation}

According to the curvature conditions and the Laplacian comparison theorem on forward complete Finsler manifolds, it satisfies the Laplacian estimation \cite{S1 2023} that
\begin{equation}
	\begin{split}
		\Delta^{\nabla u} \varphi &= \frac{\tilde{\varphi}'\Delta^{\nabla u}r}{R}+\frac{ \tilde{\varphi}'' F^{2}_{\nabla u}(\nabla^{\nabla u}r)}{R^2} \\ 
		&\geq -\frac{C_1}{R}\left[C(N,A)\sqrt{\frac{K(2R)}{C(N,A)}}~ \coth\left(R\sqrt{\frac{K(2R)}{C(N,A)}}\right)+C_0(K_{0},A)\right]-\frac{\alpha C_2}{R^2}\\
		&\geq -\frac{C_3}{R^{2}}[1+R(1+\sqrt{K(2R)})],
	\end{split}
\end{equation}
where $C_3=C_3(N,A,K_0)$ is a constant depending on $K_0$, $N$ and $A$.

Defining $Z=tH(x,t)$, we suppose that the support of function $\varphi(x)Z(x,t)$ is contained in $B_{p}(2R)$. For any fixed $T>0$, let $(x_{0},t_{0})$ be the point where $\varphi(x)Z(x,t)$ achieves its positive maximum, at which it satisfies that
\begin{equation}
	\nabla(\varphi Z)=0,   \quad \frac{\partial(\varphi Z)}{\partial_t}\geq 0  , \quad \Delta(\varphi Z)\leq 0.
	\label{4.6}
\end{equation}
(\ref{4.6}) implies that 
\begin{equation}
	\nabla Z = -\frac{\nabla \varphi}{\varphi}Z,
	\label{4.9}
\end{equation}
and 
\begin{equation}
	\Delta \varphi\cdot Z+2g_{\nabla u}(\nabla ^{\nabla u}\varphi , \nabla^{\nabla u} Z) +\varphi(\Delta-\frac{\partial}{\partial_t})Z \leq 0,
	\label{4.11}
\end{equation}
at  the maximum point $(x_0,t_0)$.
Let $\beta=\frac{a}{q}$, $a=sq^{2}$, employing Lemma \ref{lemma}, we have
\begin{equation}
	\begin{split}
		(\Delta-\frac{\partial}{\partial_t})Z &= t(\Delta-\frac{\partial}{\partial_t})H - H \\
		&\geq\left[\frac{2(1-\varepsilon)(sq+s-1)^{2}}{N-\varepsilon(N-n)s^2q^2}-2\left(\frac{1}{\varepsilon}-1\right)\right]\frac{F^{4}_{\nabla u}(\nabla ^{\nabla u}W)}{W^4}t \\
		&+ \frac{2(1-\varepsilon)}{N-\varepsilon(N-n)}\frac{1}{s^2q^2}Z^2\frac{1}{t} 
		+ \frac{4(1-\varepsilon)}{N-\varepsilon(N-n)}\frac{sq+s-1}{\varepsilon^2q^2}\frac{F^{2}_{\nabla u}(\nabla ^{\nabla u}W)}{W^2}Z\\& + \frac{2}{q}g_{\nabla u}(\nabla^{\nabla u}Z,\nabla ^{\nabla u}\log W) + cW^{-\frac{1}{q}}Z + (c-2cs)\frac{F^{2}_{\nabla u}(\nabla ^{\nabla u}W)}{W^2}W^{-\frac{1}{q}}t\\& + cW^{-\frac{1}{q}}Z - 2K(2R)\frac{F^{2}_{\nabla u}(\nabla ^{\nabla u}W)}{W^2}t - \frac{Z}{t}.
	\end{split}
	\label{4.12}
\end{equation}
It follows from the h\"{o}lder inequality that
\begin{equation}
	\begin{split}
		2K(2R)t\frac{F^{2}_{\nabla u}(\nabla ^{\nabla u}W)}{{W^2}} \leq& \frac{2(1-\varepsilon)}{N-\varepsilon(N-n)}\frac{(s-1)^2}{s^2q^2}\frac{F^{4}_{\nabla u}\left(\nabla ^{\nabla u}W\right)}{{W}^4}t\\ &+\frac{N-\varepsilon(N-n)}{2(1-\varepsilon)}\frac{s^2q^2}{(s-1)^2}K^2t,
		\label{4.13}
	\end{split}
\end{equation}
and
\begin{equation}
	\begin{split}
		(c-2cs) \frac{F^{2}_{\nabla u}(\nabla ^{\nabla u}W)}{W^2} W^{-\frac{1}{q} } t\geq -\left(\frac{(2s-1)^{2} }{4} \frac{F^{4}_{\nabla u}(\nabla ^{\nabla u} {W})}{W^4}t+c^{2}M_{1}^{2}t\right).
	\end{split}
	\label{4.14}
\end{equation}

Substituting (\ref{4.13}) and (\ref{4.14}) into (\ref{4.12}), and choosing $s>1$, $q>0$, such that $\frac{2(1-\varepsilon )}{n}\frac{s-1}{sq}\geq \frac{1}{\varepsilon }-1+\frac{(2s-1)^{2} }{8}$, one may find that
\begin{equation}
	\begin{split}
		(\Delta-\frac{\partial}{\partial_t})Z \geq&\frac{2(1-\varepsilon)}{N-\varepsilon(N-n)}\frac1{s^2q^2}Z^2\frac{1}{t}+\frac{4(1-\varepsilon)}{N-\varepsilon(N-n)}\frac{sq+s-1}{s^2q^2}\frac{F^{2}_{\nabla u}\left(\nabla^{\nabla u} W\right)}{W^2}Z \\
		&+\frac{2}{q}g_{\nabla u}\left(\nabla^{\nabla u}Z,\nabla ^{\nabla u}\log W\right)-\frac{Z}{t}- \frac{N-\varepsilon(N-n)}{2(1-\varepsilon)}\frac{s^2q^2}{(s-1)^2}K^2t-c^2M_{1}^{2}t.\\
	\end{split}
	\label{4.15}
\end{equation}
Substituting (\ref{4.15}) into (\ref{4.11}) and noticing (\ref{4.9}), we have
\begin{equation}
	\begin{split}
		&\frac{2(1-\varepsilon )}{N-\varepsilon(N-n)}\frac{1}{s^{2}q^{2}  }Z^{2}\frac{1}{t}\varphi +\frac{4(1-\varepsilon)}{N-\varepsilon(N-n)}\frac{sq+s-1}{s^2q^2}\frac{F^{2}_{\nabla u}\left(\nabla ^{\nabla u}W\right)}{W^2}Z\varphi \\
		&-\frac{2Z}{q}g_{\nabla u}(\nabla^{\nabla u}\varphi ,\frac{\nabla ^{\nabla u} W}{W} ) -\frac{Z\varphi }{t} -\frac{N-\varepsilon(N-n)}{2(1-\varepsilon)}\frac{s^2q^2}{(s-1)^2}K^{2}(2R)\varphi t \\
		&-C^2M_{1}^{2}t\varphi-\frac{C_3}{R^{2}}[1+R(1+\sqrt{K(2R)})]Z-\frac{2\alpha C_{1}^{2}  }{R^{2} } Z\leq 0.
	\end{split}
	\label{4.16}
\end{equation}
It's easy to know that
\begin{equation}
	\begin{split}
		\frac{2Z}{q}g_{\nabla u}(\nabla^{\nabla u}\varphi ,\frac{\nabla ^{\nabla u} W}{ W} ) \leq& \frac{4(1-\varepsilon)}{N-\varepsilon(N-n)}\frac{sq+s-1}{s^2q^2}\frac{F^{2}_{\nabla u}\left(\nabla W\right)}{W^2}Z\varphi\\
		&+\frac{N-\varepsilon(N-n)}{4(1-\varepsilon)}\frac{s^{2}Z }{sq+s-1} \frac{F^{2}_{\nabla u}( \nabla^{\nabla u}  \varphi)}{\varphi }.
	\end{split}
	\label{4.17}
\end{equation}
Multiplying through by $t\varphi$ at (\ref{4.16}) and utilizing (\ref{4.17}), we get that  
\begin{equation}
	\begin{split}
		&\frac{2(1-\varepsilon )}{N-\varepsilon(N-n)}\frac{1}{s^{2}q^{2}  }Z^{2}\varphi^{2} -Z\varphi -t\Bigg\{\frac{2\alpha C_{1}^{2}  }{R^{2} }\\
		&+\frac{N-\varepsilon(N-n)}{4(1-\varepsilon)}\frac{s^2}{(sq+s-1)}\frac{\alpha C^{2}_{1}  }{R^{2} }+\frac{C_3}{R^{2}}[1+R(1+\sqrt{K(2R)})] \Bigg\}Z\varphi\\
		&-t^{2}\left[\frac{N-\varepsilon(N-n)}{2(1-\varepsilon)}\frac{s^2q^2}{(s-1)^2}K^2+C^{2}M^{2}_{1} \right]\leq 0.    
	\end{split}
	\label{4.18}
\end{equation}
which implies that
\begin{equation}
	\begin{split}
		Z\varphi
		&\leq \frac{[N-\varepsilon(N-n)]s^2q^2 }{2(1-\varepsilon)}\Bigg\{1+t\Bigg(\frac{N-\varepsilon(N-n)}{4(1-\varepsilon)}\frac{s^2}{(sq+s-1)}\frac{C_{1}^{2} }{R^{2} }\\
		&\quad +\frac{C_3}{R^{2}}[1+R(1+\sqrt{K(2R)})]+\frac{2C_{1}^{2} }{R^{2} }\Bigg)\Bigg\}\\
		&\quad+t\sqrt{\frac{[N-\varepsilon(N-n)]s^2q^2 }{2(1-\varepsilon)}}\sqrt{\frac{[N-\varepsilon(N-n)]s^2q^2 }{2(1-\varepsilon)(s-1)^2}AK^2+c^2M^2_{1}}, 
	\end{split}
\end{equation}
is satisfied at $(x_{0},t_{0})$.

Clearly, corresponding to the assumption of $(x_0,t_0)$, one can deduce that
\begin{equation}
	\begin{split}
		&\sup_{x\in B_{p}} Z(x,T) \leq Z(x_0,t_0)\varphi(x_0) \\
		\leq& \frac{[N-\varepsilon(N-n)]s^2q^2 }{2(1-\varepsilon)}\Bigg\{1+T\Bigg(\frac{N-\varepsilon(N-n)}{4(1-\varepsilon)}\frac{s^2}{(sq+s-1)}\frac{C_{1}^{2} }{R^{2} }+\frac{2C_{1}^{2} }{R^{2} }\\
		&\quad\quad +\frac{C_3}{R^{2}}[1+R(1+\sqrt{K(2R)})]\Bigg)\Bigg\}\\
		&+T\sqrt{\frac{[N-\varepsilon(N-n)]s^2q^2 }{2(1-\varepsilon)}}\sqrt{\frac{[N-\varepsilon(N-n)]s^2q^2 }{2(1-\varepsilon)(s-1)^2}AK(2R)^2+c^2M^2_{1}}.
	\end{split}
\end{equation}
Recall a result that, if $a_{1}, a_{2}, a_{3} \geq 0$ then $\sqrt{a_{1}}\sqrt{a_{2}+a_{3}} \leq \sqrt{a_{1}a_{2}}+\sqrt{a_{1}a_{3}}$. Thus,
\begin{equation}
	\begin{split}
		&\frac{F^2_{\nabla u}(\nabla u)}{u^2}+sc(1-u) -s\frac{u_{t}}{u}\\
		\leq &\frac{N-\varepsilon(N-n)}{2(1-\varepsilon)}\frac{s^{2}}{t} +\frac{N-\varepsilon (N-n)}{2(1-\varepsilon )(s-1)} s^{2}K(2R)+\frac{sc}{q} \sqrt{\frac{N-\varepsilon (N-n)}{2(1-\varepsilon )}}M_{1} \\
		&+\frac{[N-\varepsilon (N-n)]s^{2}}{2(1-\varepsilon )R^{2}} \Bigg\{C_3[1+R(1+\sqrt{K(2R)})]+\frac{[N-\varepsilon(N-n)]s^2}{4(1-\varepsilon)(s-1)}C^{2}_{1}+2C_{1}^2\Bigg\},	
	\end{split}
\end{equation} 
on $B(p,2R)$, which is the desired inequality.  
\end{proof}

\begin{remark}
	Theorem  \ref{Theorem 1.2} follows from Theorem \ref{theorem of noncompact} directly.
\end{remark}



\end{document}